\documentclass[a4paper, 10pt, footinclude=false]{article}
\usepackage{geometry}
\usepackage{amsthm, amsmath, amssymb, amsfonts}
\usepackage[utf8]{inputenc}
\usepackage[english]{babel}
\usepackage{url}
\usepackage{mathtools}

\usepackage{bm}
\usepackage{tikz}
\usepackage{pgfplots}
\usepackage{pgfplotstable}
\usepackage{booktabs}
\usepackage{stackengine}
\usepackage{algorithm}
\usepackage[noend]{algpseudocode}
\usepackage{layouts}

\numberwithin{figure}{section}
\numberwithin{table}{section}
\numberwithin{equation}{section}

\newenvironment{abstr}[1]{ \vspace{.05in}\footnotesize
	\parindent .2in
	{\upshape\bfseries #1. }\ignorespaces}{\par\vspace{.1in}}

\DeclareMathAlphabet\mathbfcal{OMS}{cmsy}{b}{n}

\newtheorem{theorem}{Theorem}[section]
\newtheorem{lemma}[theorem]{Lemma}

\newtheorem{assumption}[theorem]{Assumption}

\theoremstyle{definition}

\newtheorem{remark}[theorem]{Remark}
\newtheorem{example}[theorem]{Example}

\newcommand{\RR}{\mathbb{R}}
\newcommand{\NN}{\mathcal{N}}

\newcommand{\II}{\mathcal{I}}

\newcommand{\hatVV}{\mathbf{\hat V}}
\newcommand{\VV}{\mathbf{V}}

\newcommand{\uu}{\mathbf{u}}
\newcommand{\vv}{\mathbf{v}}
\newcommand{\ww}{\mathbf{w}}

\newcommand{\KK}{\mathbf{K}}
\newcommand{\ff}{\mathbf{f}}
\newcommand{\zz}{\mathbf{z}}
\renewcommand{\gg}{\mathbf{g}}
\newcommand{\LL}{\mathbf{L}}
\newcommand{\MM}{\mathbf{M}}
\newcommand{\PP}{\mathbf{P}}
\newcommand{\BB}{\mathbf{B}}

\newcommand{\III}{\mathbfcal{I}}

\newcommand{\spanop}{\operatorname{span}}

\newcommand{\vol}{\operatorname{vol}}

\newcommand{\norm}[1]{{\left\vert\kern-0.25ex\left\vert\kern-0.25ex\left\vert #1 
    \right\vert\kern-0.25ex\right\vert\kern-0.25ex\right\vert}}

\usepackage{todonotes}

\allowdisplaybreaks[4]

\begin{document}
	
	\title{Iterative solution of spatial network models by subspace decomposition}
	\author{
	M. G\"{o}rtz\textsuperscript{1}, F. Hellman\textsuperscript{2}, 
	A. M\aa lqvist\textsuperscript{2}}
	
	\date{}
	
\maketitle
\footnotetext[1]{Fraunhofer-Chalmers Centre, Chalmers Science Park, 412 88 G\"{o}teborg, Sweden} \
\footnotetext[2]{Department of Mathematical Sciences, Chalmers University of Technology and University of Gothenburg, 412 96 G\"oteborg, Sweden}

\begin{abstract}
We present and analyze a preconditioned conjugate gradient method (PCG) for solving spatial network problems. Primarily, we consider diffusion and structural mechanics simulations for fiber based materials, but the methodology can be applied to a wide range of models, fulfilling a set of abstract assumptions. The proposed method builds on a classical subspace decomposition into a coarse subspace, realized as the restriction of a finite element space to the nodes of the spatial network, and localized subspaces with support on mesh stars. The main contribution of this work is the convergence analysis of the proposed method. The analysis translates results from finite element theory, including interpolation bounds, to the spatial network setting. A convergence rate of the PCG algorithm, only depending on global bounds of the operator and homogeneity, connectivity and locality constants of the network, is established. The theoretical results are confirmed by several numerical experiments.
\end{abstract}
\maketitle

\section{Introduction}

\begin{figure}
 \centering
 \includegraphics{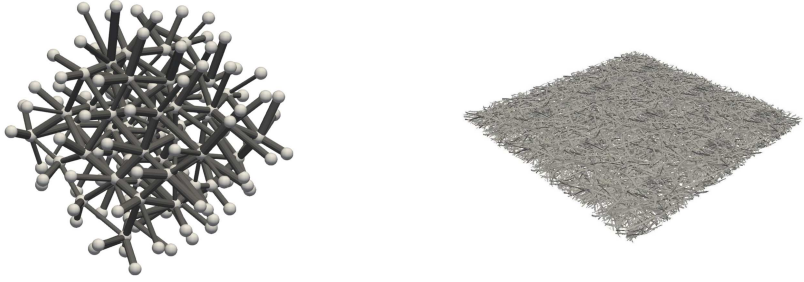}
  \caption{Illustration of a subgraph in a pore network model (left) and a fiber network model (millimeter scale) of paper (right).}
  \label{fig:porenetworkandpaper}
\end{figure}

Many phenomena in science and engineering, modelled by partial differential equations (PDE), are challenging to simulate accurately due to their vast complexity. Sometimes it is beneficial to simplify, still maintaining the main features of the full model, by introducing a spatial network model. 
In porous media flow modelling, see Figure~\ref{fig:porenetworkandpaper} (left), the flow in the exact pore geometry can be approximated by a simplified network model of edges (throats) and nodes (pore cavities), see \cite{Blunt,Gjennestad,Huang}. This technique reduces the computational complexity and allows for simulation using larger computational domains. A fiber based material such as paper, see Figure~\ref{fig:porenetworkandpaper} (right), can also be modelled using a spatial network model where the three dimensional hollow cylindrical fibers are modelled as one dimensional objects, represented by edges connected with nodes, and the contact regions as additional nodes, see \cite{Heyden,Iliev,BIT}. Even though much of the complexity is eliminated in this way, it is still often challenging to solve the resulting model efficiently.

In this paper, we consider spatial network models that arise from applications modelled by linear elliptic partial differential equations such as heat conduction and structural mechanics. The corresponding system matrices are weighted graph Laplacians in the scalar case and related to the weighted graph Laplacians in each coordinate direction in the vector valued case. Successful numerical algorithms for solving the resulting sparse linear systems use parallelization and discretizations on a range of scales. This is exploited in iterative methods such as geometric and algebraic multigrid \cite{multigrid,AMG,XZ17} and domain decomposition \cite{Saad, PB96}. In a PDE setting with a simple geometry it is natural to introduce nested discretizations on different levels. In a purely algebraic setting it is less obvious but there are many coarsening techniques available that uses connectivity information from the system matrix, see for instance \cite{AMG,XZ17}. There are geometry based algebraic multigrid methods that uses some geometrical information to construct levels of coarsening \cite{XZ17}. 
This setting is most similar to the model we study.
    
We consider a spatial network defined by a graph $\mathcal{G}=\{\mathcal{N},\mathcal{E}\}$ of nodes and edges. The nodes and edges are contained in a bounded domain $\Omega$. The model problem is defined by a symmetric linear operator $\KK$ and a right hand side $\ff$. The operator act on a linear space $\VV$ of functions defined on $\mathcal{N}$. We want to solve a stationary equation of the form: find $\uu\in\VV$ such that for all $\vv\in \VV$,
$$
(\KK \uu,\vv)=(\ff,\vv).
$$
We introduce a coarse finite element discretization of $\Omega$ and let $\VV_H\subset\VV$ be its restriction to the network.  The full solution space $\VV$ is then decomposed into $\VV_0=\VV_H$ and a set of overlapping local subspaces $\{\VV_j\}_{j=1}^m$. Given the subspace decomposition we formulate an additive Schwarz preconditioner for the model problem and apply the conjugate gradient method to the preconditioned system. The main theoretical contribution of this work is the convergence proof of the proposed scheme, with a rate that only depends on global bounds on $\KK$ and the homogeneity, connectivity and locality of the network on a scale $R_0$, under the additional assumption that $H\geq R_0$. These dependencies are traced explicitly in the analysis and reveals the impact of the $R_0$-scale homogeneity and connectivity on the convergence of the iterative solver. To prove the main result we translate results from finite element theory to the network setting including the construction of an interpolation operator, the proof of interpolation bounds, and the proof of a product rule type bound. The interpolation bound in turn builds on Poincar\'{e} and Friedrichs inequalities on subgraphs, see \cite{ChungSurvey, ChungBook} and the earlier contributions \cite{cheeger,Fiedler}.  By this we establish spectral bounds of the preconditioned operator. Given the  spectral bounds we apply classical additive Schwarz theory to prove the convergence, see \cite{KY16, Xu90}. In the first numerical experiment we consider a randomly generated model of a fiber based material, defined by fiber length and density. We investigate how the homogeneity and connectivity constants depends on the parameters in the model and the size of the domain. Then we present a scalar (heat conductivity) and a vector (structural mechanical model) valued numerical example which shows the efficiency of the proposed method for models of fiber based materials.

The paper is organized as follows. Section 2 is devoted to preliminary notation and Section 3 to problem formulation and assumptions on the operator $\KK$ and the network. In Section 4 we introduce a preconditioned conjugate gradient method and present the main convergence result.  Section 5 is devoted to the proof of the spectral bound of the preconditioned operator which is the main technical result needed for the convergence proof.
Finally, in Section 6 we present numerical examples. 

\section{Preliminary notation}

We start by introducing notation for the spatial network itself and for operators defined on the network. We also introduce finite element spaces that are restricted to the network and used in the subspace decomposition.

\subsection{Spatial network}
\begin{figure} [t]
 \centering
 \includegraphics{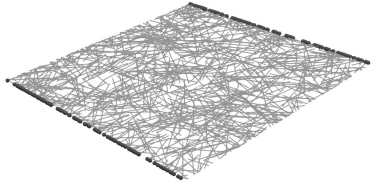}
  \caption{Dirichlet nodes on a subset $\Gamma$ of the boundary of the domain $\Omega$.}
  \label{fig:gamma}
\end{figure}

We consider a connected spatial network represented as the graph
$\mathcal{G} = (\mathcal{N}, \mathcal{E})$. The node set
$\mathcal{N}$ is a finite set of points in $\RR^d$ and the edge set
$$\mathcal{E} = \{ \{x, y\} \,:\, \text{an edge connects } x
\text{ and } y \in \mathcal{N} \}$$ consists of unordered pairs of the
endpoint vertices of each edge. The notation $x \sim y$ expresses that
$\{x, y\}$ is an edge in $\mathcal{E}$. If $x \sim y$, then we say
that the nodes $x$ and $y$ are adjacent. If two nodes $x$ and $y$ are adjacent
then $| x - y |$ is the length of the edge connecting them, with
$| \cdot |$ being the Euclidean norm. The network is embedded in a
spatial domain $\Omega$ ($\mathcal{N} \subset \Omega$), which we for
simplicity assume to be the closed hyper-rectangle
$$\Omega = [0, l_1]\times[0,l_2]\times\dots\times [0,l_d] \subset \mathbb{R}^d.$$ In Remarks \ref{rem:poly1}
we discuss how the methodology 
can be modified for general poly\-gonal/poly\-hedral
domains.  We let $\Gamma \subset \partial \Omega$ be the part of the
boundary of $\Omega$ on which homogeneous Dirichlet boundary
conditions are applied, see Figure~\ref{fig:gamma}.

We let $\hat V$ be the space of real-valued functions defined on the
node set $\NN$ and let
$$V = \{ v \in \hat V\,:\, v(x) = 0, x \in \Gamma\}$$
denote the function set satisfying the boundary conditions. For two
functions $u, v \in \hat V$, we define the node-wise product
$uv \in \hat V$ in the natural way i.e. $(uv)(x) = u(x)v(x)$ for
$x \in \mathcal{N}$. When $\omega \subset \Omega$, we
define $\mathcal{N}(\omega) = \mathcal{N} \cap \omega$.

\subsection{Operators and norms}
For $u,v \in \hat V$ we define the inner product
$$(u, v) = (u, v)_\Omega \quad\text{where}\quad (u, v)_\omega = \sum_{x \in \mathcal{N}(\omega)} u(x) v(x).$$
For every node in the network, $x \in \mathcal{N}$, let
$M_x: \hat V \to \hat V$ be the diagonal linear operator defined by
\begin{equation}
  (M_x v, v) = \frac{1}{2}\sum_{x \sim y} |x - y| v(x)^2.
\end{equation}
By $\sum_{x \sim y}$ we mean that the sum goes over all nodes $y$
adjacent to the given $x$.  For subdomains $\omega \subset \Omega$ we
introduce the shorthand notation
$M_\omega = \sum_{x \in \mathcal{N}(\omega)} M_x$, and $M =
M_\Omega$. Note that $M:\hat V\rightarrow\hat V$ but, if we restrict
the domain, also $M:V\rightarrow V$ holds. We define the norm
$|v|_M = (Mv, v)^{1/2}$ and the semi-norms
$|v|_{M,\omega} = (M_\omega v, v)^{1/2}$.

Next we define the reciprocal edge length weighted graph Laplacian. Let
$L_x : \hat V \to \hat V$ be defined by
\begin{equation}
  (L_x v, v) = \frac{1}{2}\sum_{x\sim y} \frac{(v(x) - v(y))^2}{|x - y|}.
\end{equation}
Again we introduce the shorthand notation for subdomains
$\omega \subset \Omega$, $L_\omega = \sum_{x \in \mathcal{N}(\omega)} L_x$, and
$L = L_\Omega$.  Note that
$(L_\omega u)(x)$ is generally nonzero for vertices $x$ outside
$\omega$ that are adjacent to a node in $\omega$. The weighted
Laplacian $L$ is symmetric and positive semi-definite and the kernel of $L$ contains the 
constant functions in $\hat V$. Since we
assume that the network is connected, the kernel of $L$ has 
dimension one. The operator $L$ defines the following semi-norms,
$|v|_L = (Lv, v)^{1/2}$ and $|v|_{L,\omega} = (L_\omega
v, v)^{1/2}$.
The notation $V(\omega)$ will be used to denote the space of functions
that are zero for nodes outside $\omega$.

To handle both scalar and vector valued functions we introduce vector-valued functions of $n$ components, where for example $n=1$ or $n=d$.
We introduce the product space
\begin{equation*}
\VV = V^n = V \times \cdots \times V \qquad \text{(}n\text{ times)}
\end{equation*}
as the admissible function space for the unknown and the full space
$\hatVV = {\hat V}^n$ (so that $\VV \subset \hatVV$). The components
of $\VV$ need not be identical if different boundary conditions are
applied to the components, however, for simplicity we assume that all
components are identical. A function $\vv \in \hatVV$ consists of the
components $\vv = [v_1, v_2, \ldots, v_n]$ subscripted by
their index. We introduce $\LL_x : \hat\VV \to \hat\VV$
as $L_x$ applied componentwise, i.e.
$$\LL_x\vv = [L_xv_{1}, \ldots, L_xv_{n}].$$ For
$\omega \in \Omega$, we let
$\LL_\omega = \sum_{x \in \mathcal{N}(\omega)} \LL_x$ and 
$\LL = \LL_\Omega$. We also extend the notation for the inner product
to the product space by
$$(\uu, \vv) = (u_{1}, v_{1}) + \cdots + (u_{n}, v_{n})$$ and
introduce the semi-norms
$$|\vv|_{\LL,\omega} = (\LL_\omega \vv, \vv)^{1/2} =
\left(|v_1|_{L,\omega}^2 + \cdots + |v_n|_{L,\omega}^2\right)^{1/2}$$
and shorthand $|\vv|_{\LL} = |\vv|_{\LL,\Omega}$. The analogue
notation for $\MM$ will also be used.

\subsection{Finite element mesh and function space}\label{sec:mesh}
The preconditioner we propose for the iterative solver uses a coarser
scale representations of the full network. To construct this
representation we first introduce a family of finite element meshes on
the spatial domain
$\Omega$. The construction we consider has been used previously to handle 
non-nested meshes in the localized orthogonal decomposition 
method, see \cite{bookLOD,MaPe15}.
To be able to emphasize the main message of the paper, we choose a
simple mesh of hypercubes (squares for $d = 2$, and cubes for $d = 3$,
etc.).

We introduce 
boxes $B_R(x) \subset \Omega$ centered at $x = (x_1, \ldots, x_d)$
with side length $2R$ as follows. Let
$$B_R(x) = [x_1 - R, x_1 + R) \times \cdots \times [x_d - R, x_d
+ R),$$ but replace $[x_i - R, x_i + R)$ with $[x_i - R, x_i + R]$
for any $i$ for which $x_i + R = l_i$. In this way,
$B_R(x)$ only include its upper boundary of any dimension if it intersects the boundary of $\Omega$. We let
$\mathcal{T}_H$ be a family of partitions (meshes) of $\Omega$ into
hypercubes (elements) of side length $H$:
\begin{equation*}
  \mathcal{T}_{H} =
  \{B_{H/2}(x) \,:\, x = (x_1, \ldots, x_d) \in \Omega \text{ and } H^{-1}x_i + 1/2 \in \mathbb{Z} \text{ for } i = 1, \ldots, d\}.
\end{equation*}

\begin{figure}[t]
  \centering
  \includegraphics{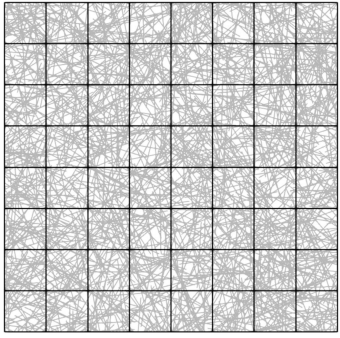}
  \caption{The finite element mesh $\mathcal{T}_H$ with a network in the background.}
  \label{fig:meshondomain}
\end{figure}

For simplicity, we assume the side lengths of the hyper-rectangle $l_1,l_2,\dots,l_d$ are integer multiples of a mesh parameter $H$ so that $\mathcal{T}_H$ exactly covers $\Omega$. By the construction of $B_{R}(x)$,  $\mathcal{T}_H$ is a true partition so that
each point in $\Omega$ is included in exactly one element, see Figure~\ref{fig:meshondomain}.

\begin{figure}[b]
 \centering
 \includegraphics{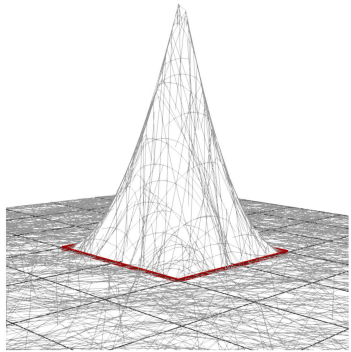}
 \caption{The basis function $\varphi_j$ with corresponding patch $U(y_j)$ in red.} \label{fig:UandPhi}
\end{figure}

To name patches of elements in a mesh $\mathcal{T}_H$, we introduce
the notation $U$. For $\omega \subset \Omega$ let
$$U(\omega) = \{ x \in \Omega \,:\, \exists T \in \mathcal{T}_H \,:\, x \in T, \ \overline T \cap \overline \omega
\ne \emptyset\}.$$ The patch $U(T)$ contains
the points in $T$ and its adjacent elements.  Recursively, we define larger patches $U_j(\omega) = U_{j-1}(U(\omega))$ with
$U_1 = U$. By this construction, $U_2(T) = U(U(T))$ contains the points in all elements in the
element patch for $T$ with ``radius'' 2. Moreover, the definition of $U$ is extended to accept arguments that are individual points, by $U(x) = U(\{x\})$ when $x \in \Omega$, see Figure~\ref{fig:UandPhi}.

We now introduce a first-order finite element space on the mesh
$\mathcal{T}_H$. Let $\mathcal{\hat Q}_H$ be the space of continuous
functions defined on $\Omega$ whose restriction to
$T \in \mathcal{T}_H$ is a linear combination of the polynomials
$z = (z_1, \ldots, z_d) \mapsto z^{\alpha}$ for multi-index $\alpha$ with $|\alpha| \le 1$. For $d = 2$,
this is the space of bilinear functions on $T$. Further, let
$\mathcal{Q}_{H} = \{ q \in \mathcal{\hat Q}_H \,:\, q(x) = 0, x \in \Gamma \}$.

We are interested in the restriction of these functions to the
network nodes. Let $\hat V_H$ be the functions in
$\mathcal{\hat Q}_H$ restricted to $\mathcal{N}$, with $V_{H}$ similarly defined. We further define the function space $\VV_H$ of vector-valued functions by
$$
\VV_H=V_H\times\cdots\times V_H  \qquad \text{(}n\text{ times)}.
$$

From now on we consider a fix $H$ and omit subscript $H$ in new notation. 
Let $m$ be the dimension of $\mathcal{\hat Q}_H$, and denote the Lagrange finite element basis functions by $\varphi_1, \ldots, \varphi_m \in \mathcal{\hat Q}_{H}$,
and their restrictions to the network nodes again denoted by $\varphi_1, \ldots, \varphi_m \in \hat V_{H}$, see Figure~\ref{fig:UandPhi}. The nodes of the mesh $\mathcal{T}_H$ are denoted  $y_1, \ldots, y_m$. Note that we do not require that the mesh nodes coincide with network nodes in $\mathcal{N}$.

We assume that the mesh is aligned with $\Gamma$ in the sense that $\Gamma$
can be expressed as a union of element faces.  Without loss of
generality we order the basis function such that the $m_0 < m$
first basis functions $\varphi_1, \ldots, \varphi_{m_0}$ span
$\mathcal{Q}_{H}$.

\begin{remark}\label{rem:poly1}
If the domain $\Omega$ is a general polygonal/polyhedral domain we instead use a finite element space based on continuous piecewise affine functions on a shape regular mesh of $d$-simplices (triangles for $d=2$ and tetrahedra for $d=3$). Every node in the network is assigned to exactly one element. If a point belongs to the closure of several elements we assign it to the one with the lowest element number. Once the elements are constructed and the nodes are distributed to the elements, element patches $U$ can again be defined as above.     
\end{remark}

\section{Problem formulation and assumptions}

In this section we introduce the model problem. We also state assumptions on the linear operator $\KK$, defining the interactions over the network, and the network itself needed to prove our main convergence result. 

The model problem is expressed in terms of a symmetric and positive semi-definite linear operator
$\KK : \hat\VV \to \hat\VV$ and a right hand side function $\ff \in \hat \VV$: find
$\uu \in \VV$ such that for all $\vv \in \VV$,
\begin{equation}\label{eq:weakform}
  (\KK \uu, \vv) = (\ff, \vv).
\end{equation}

Since $\uu \in \VV$, the solution is zero
in the non-empty set of nodes $\mathcal{N}(\Gamma)$ where
$\Gamma \subset \partial \Omega$ is the Dirichlet boundary. If the
sought solution is equal to some given non-zero function
$\gg(x)$ for $x \in \mathcal{N}(\Gamma)$, we extend
$\gg$ to all nodes, let the sought solution be $\hat \uu=\uu+\gg$, and solve for
$\uu \in \VV$ with modified right hand side instead. Non-zero flux boundary conditions can be added in the right hand side of equation (\ref{eq:weakform}).

\subsection{Assumptions on $\KK$}
We restrict our attention to operators with the following properties:
\begin{assumption}\label{ass:K} The operator $\KK$
  \begin{enumerate}
  \item is bounded and coercive on $\VV$ with respect to $\LL$, i.e.\
    there are constants $\alpha > 0$ and $\beta < \infty$ such that
    \begin{equation}\label{eq:equivKL}
      \alpha (\LL \vv, \vv) \le (\KK \vv, \vv) \le \beta (\LL \vv, \vv)
    \end{equation}
    for all $\vv \in \VV$, 
\item is symmetric, $(\KK\vv,\ww)=(\KK\ww,\vv)$, and
\item  admits a unique solution to equation (\ref{eq:weakform}).
  \end{enumerate}
\end{assumption}
Assumption \ref{ass:K}.3 implicitly puts a constraint on $\VV$. A sufficient portion of the boundary has to be fixed in order for it to hold. In the simple case of $n=1$ and $\KK=L$ it is enough that $\Gamma$ contains one node. In the vector valued case the kernel is larger and more nodes need to be fixed. We note that the bilinear form $(\KK \cdot,\cdot)$ is an inner product on $\VV$. 

\begin{example}[Heat conductivity]
  This is a scalar example with $n=1$ and we thus drop the bold-face of $\KK$ and $\uu$. Here $u$ is the scalar temperature distribution in the nodes of
  the network. Then we can define $K=\sum_{x \in \mathcal{N}}K_x$ with
  \begin{equation}\label{eq:scalarmodel}
    (K_xv,v) = \frac{1}{2}\sum_{x \sim y} \gamma_{xy} \frac{(v(x)-v(y))^2}{|x - y|},
  \end{equation}
  where $0 < \gamma_{xy} < \infty$ is a heat conductivity material
  parameter for the connecting edges. Assumption~\ref{ass:K} is
  satisfied with $\alpha = \min_{x \sim y} \gamma_{xy}$ and $\beta = \max_{x \sim y} \gamma_{xy}$. The operator $K$ is clearly symmetric by construction. By putting at
  least one of the nodes at the Dirichlet boundary $\Gamma$, we remove
  the constants from $V$ and the kernel of $L$ restricted to $V$
  contains only zero and we have a unique solution. The right hand
  side $f$ (given for each node) is the external heat source.
\end{example}

\begin{example}[Structural problem for fiber based material] \label{ex:fiber}
  This is a vector-valued example with $d=n=3$.
  Let the network edges represent a composite material made up of fibers. The forces resulting from a specific displacement of the network can be expressed as the model problem in \eqref{eq:weakform}
  where $\uu(x)$ represents the displacement of the network node $x$ and $\textbf{f}(x)$ the resulting directional forces in that node. The model presented in this example is described in detail in \cite{BIT}, and is a linearization of Hooke's law and Euler--Bernoulli beam theory. 
  
  The model $\KK$ can be decomposed into two parts as follows,
  $$\KK = \KK^{(\text{T})} + \KK^{(\text{B})} = \sum_{x\in\mathcal{N}} \KK^{(\text{T})}_x + \KK^{(\text{B})}_x,$$
  where $\KK_x^{(\text{T})} $ represents the tensile part (Hooke) and $\KK_x^{(\text{B})}$ the bending stiffness contributions (Euler--Bernoulli). 
  
  The tensile stiffness addition can be defined by the following bilinear form,
  \begin{equation} \label{eq:tensile}
    (\KK_x^{(\text{T})}\vv,\vv) = \frac{1}{2}\sum_{x \sim y} \gamma_{xy}  \frac{((\vv(x)-\vv(y)) \cdot \partial_{xy})^2}{|x - y|}, 
  \end{equation}
  where $\gamma_{xy}$ is the tensile stiffness, $\partial_{xy}=|x-y|^{-1}(x-y)$ is the directional vector for the corresponding edge.
      
  Similar to the tensile stiffness addition, the bending stiffness contribution can be defined by the following bilinear form,
  \begin{equation} \label{eq:bend}
  \begin{split}
    & (\KK^{(\text{B})}_x\vv,\vv) = \\
    & \sum_{\substack{\ k \in \{1, \ 2 \} \\x\sim y \wedge x\sim z \\ y \not=z}} \gamma_{xyz}^{(k)} \frac{|x-y|+|x-z|}{2}  \left(\frac{(\vv(y)-\vv(x)) \cdot \eta_{xyz}^{y,(k)}}{|x-y|}+\frac{(\vv(z)-\vv(x)) \cdot \eta_{xyz}^{z,(k)}}{|x-z|}\right)^2,
  \end{split}
  \end{equation}
   where $\gamma_{xyz}^{(k)}$ is a constant, $\eta_{xyz}^{y,(1)}=\eta_{xyz}^{z,(1)}$ is an orthogonal unit vector to $\partial_{xy}$ and $\partial_{xz}$, and $\eta_{xyz}^{w,(2)} = \partial_{xw}\times \eta_{xyz}^{w,(1)}$ for $w = y$ and $z$. 
   
   This model satisfies the assumptions in Assumption~\ref{ass:K} conditionally. The constant $\beta$ in the first assumption is dependent on the constants and difference in edge lengths between two edges that share the same node. Whether there exists an $\alpha$, and the value of such a constant, depends on the constants and the geometry of the network. The operator is symmetric, and for invertibility we need at least three nodes in $\Gamma$ that span a plane.
\end{example}

\subsection{Assumptions on the network}

The network must resemble a homogeneous material on coarse scales, so we next assume
homogeneity, connectivity, locality and boundary density of the
network scales greater than $R_0$.
\begin{assumption}[Network assumptions]
  \label{ass:network}
  There is a length-scale $R_0$, a uniformity constant $\sigma$, and a
  density $\rho$, so that
  \begin{enumerate}
  \item (homogeneity) for all $R \ge R_0$ and $x \in \Omega$, it holds
    that
    \begin{equation*}
      \rho \le (2R)^{-d}|1|^2_{M,B_R(x)} \le \sigma \rho,
    \end{equation*}
    
  \item (connectivity) for all $R \ge R_0$ and $x \in \Omega$, there is
    a connected subgraph
    $\mathcal{\bar{G}} = (\mathcal{\bar{N}}, \mathcal{\bar{E}})$ of
    $\mathcal{G}$, that contains
    \begin{enumerate}
    \item all edges with one or both endpoint in $B_R(x)$,
    \item only edges with endpoints contained in $B_{R+R_0}(x)$,
    \end{enumerate}
  \item (locality) the edge length $|x - y| < R_0$ for all edges $\{x, y\} \in \mathcal{E}$,
  \item (boundary density) for any $y \in \Gamma$, there is an
    $x \in \mathcal{N}(\Gamma)$ such that
    $|x - y| < R_0$.
  \end{enumerate}
\end{assumption}

The homogeneity assumption requires that the density of the network
(in terms of the $M$-norm mass $|1|_{M,\omega}^2$) should be
comparable at scale $R_0$ over the whole domain. The
connectivity assumption requires that the network is locally connected
at scale $R_0$, while the edge length assumption requires that the
network is not connected over larger distances than $R_0$.  Finally,
the boundary density requires that the boundary nodes (the nodes in
$\mathcal{N}(\Gamma)$) occur frequently on scale $R_0$ in $\Gamma$.

\subsubsection{Friedrichs and Poincaré inequalities}
The connectivity assumption rules out networks that are
insufficiently locally connected. A graph defined by the edges of a
Peano curve, for instance, does not satisfy the connectivity
assumption on any reasonable small scale $R_0$, while a grid does. The
assumption can be used to show the following Friedrichs and Poincaré
inequalities that are necessary to establish interpolation bounds
in Lemma~\ref{lem:interpolation}.

\begin{lemma}[Friedrichs and Poincaré inequalities]
  \label{lem:poincare}
  If Assumption~\ref{ass:network} holds, then there is a
  $\mu < \infty$ such that for all $R \ge R_0$ and $x \in \Omega$ for
  which
  \begin{itemize}
  \item (Friedrichs) $B_{R}(x)$ contains boundary nodes, it holds that
    \begin{equation*}
      |v|_{M,B_R(x)} \le \mu R |v|_{L,B_{R + R_0}(x)},
    \end{equation*}
    for all $v \in V$,
  \item (Poincaré) $B_{R}(x)$ may or may not contain boundary nodes, it holds that
    \begin{equation*}
      |v - c|_{M,B_R(x)} \le \mu R |v|_{L,B_{R + R_0}(x)},
    \end{equation*}
    for some constant function $c = c(R, x, v)$, for all $v \in \hat V$.
  \end{itemize}
\end{lemma}
Note that the $L$-norm in the right hand side of the inequalities are
taken on a slightly larger box than the $M$-norm. This is necessary
since nodes in the smaller box could connect with other nodes within
the box only through edges connected to nodes outside the box. Thus,
the function values ($M$-norm) cannot generally be bounded by
differences of the function values ($L$-norm) only within the box
itself.
\begin{proof}[Proof of Lemma~\ref{lem:poincare}.]
  Consider an $x \in \Omega$ and $R \ge R_0$. Let $\bar{M}$ and
  $\bar{L}$ be defined for $\mathcal{\bar{G}}$ (a subgraph from
  Assumption~\ref{ass:network}.2) as $M$ and $L$ are defined for
  $\mathcal{G}$. If $B_R(x)$ contains a boundary node, then so does
  $\mathcal{\bar{N}}$, and we prove the Friedrich inequality in the
  Dirichlet case below. The Poincaré inequality is proven regardless
  of whether $B_R(x)$ contains boundary nodes or not in the Neumann
  case below.

  \begin{itemize}
  \item (Dirichlet, if $\mathcal{\bar{N}}$ contains boundary nodes).
    Denote by $\lambda_1 \le \lambda_2 \le \cdots$ the eigenvalues of
    $\bar{L}u = \lambda \bar{M}u$ for $u \in V(\mathcal{\bar{N}})$.
    Since the subgraph is connected and there are prescribed values
    for $u$ in the boundary nodes, the first eigenvalue
    $\lambda_1 > 0$. Then we have by the min-max theorem
    \begin{equation}
      \label{eq:minmax_dirichlet}
      \lambda_1 = \min_{\substack{v \ne 0}} \frac{(\bar L v, v)}{(\bar M v, v)}.
    \end{equation}
    With the assumption on the edges included and not included in
    $\mathcal{\bar G}$, we get
    \begin{equation}
      \label{eq:eigenvalue_dirichlet}
      |v|_{M,B_R(x)}^2 \le (\bar M v, v) \le \lambda_1^{-1} (\bar L v, v) \le \lambda_1^{-1} |v|_{L,B_{R+R_0}(x)}^2
    \end{equation}
    and the Friedrichs inequality is satisfied with constant
    $\mu(x, R) = R^{-1}\lambda_1^{-1/2}$.
    
  \item (Neumann) Denote by $\lambda_1 \le \lambda_2 \le \cdots$ the
    eigenvalues of $\bar{L}u = \lambda \bar{M}u$ for
    $u \in \hat V(\mathcal{\bar{N}})$, i.e.\ $u$ is free from boundary
    conditions.  Since the subgraph is connected and there are no
    prescribed boundary values of $u$ in this case, the first
    eigenvalue $\lambda_1 = 0$ (corresponding to constant
    eigenvectors) while $\lambda_2 > 0$.  Then
    \begin{equation}
      \label{eq:minmax_neumann}
      \lambda_2 = \min_{\substack{v \ne 0\\ (\bar{M}v, 1) = 0}} \frac{(\bar L v, v)}{(\bar M v, v)}.
    \end{equation}
    With $c$ as the $\bar M$-orthogonal projection of $v$ onto the
    constant functions we get
    \begin{equation}
      \label{eq:eigenvalue_neumann}
      \begin{split}
      |v - c|_{M,B_R(x)}^2 &\le (\bar M (v - c), v - c) \\
      &\le \lambda_2^{-1} (\bar L v, v) \le \lambda_2^{-1} |v|_{L,B_{R+R_0}(x)}^2
      \end{split}
    \end{equation}
    and the Poincaré inequality is satisfied with constant
    $\mu(x, R) = R^{-1}\lambda_2^{-1/2}$.
  \end{itemize}
  Since there are only finitely many subgraphs $\mathcal{\bar{G}}$
  over all possible values of $x \in \Omega$ and $R \ge R_0$, we can set
  $\mu = \max_{x,R} \mu(x,R) < \infty$.
\end{proof}

The constant $\mu$ should be small to bound the interpolation constant
in Lemma~\ref{lem:interpolation}.  The assumptions in
Assumption~\ref{ass:network} only guarantee that such a constant
exists, but not that it is small. As is illustrated in the proof, they can be
computed from eigenvalues of a Laplacian on the subgraphs
$\mathcal{\bar{G}}$. In the first part of Section~\ref{sec:numexp} we
show numerically that the constant $\mu$ appearing in the Poincaré
inequality is bounded and of moderate size for the class of randomly
generated networks used in the numerical experiments.

Next, we show that $\mu$ can be bounded by the isoperimetric constant
if there are subgraphs $\mathcal{\bar{G}}$ that satisfy a
$d$-dimensional isoperimetric inequality.

\subsubsection{Isoperimetric dimension and constant}
We introduce additional graph notation. Let $\bar d_x$ be the
degree (number of connected edges) of a node $x$ in a subgraph
$\mathcal{\bar{G}}$, and for any node subset
$X \subset \mathcal{\bar{N}}$, let
\begin{equation*}
  \vol(X) = \sum_{x \in X} \bar d_x.
\end{equation*}
Further, let $\mathcal{\bar{E}}(X, X') \subset \mathcal{\bar{E}}$
denote the set of subgraph edges with one endpoint in $X$ and the other in
$X'$. 

\begin{lemma}[Connectivity assumption by an isoperimetric inequality]
  \label{lem:isoperimetric}
  If Assumption~\ref{ass:network} holds, and
  there are constants $\nu_1$ and $\nu_2$ such that
  \begin{enumerate}
  \item[2\textquotesingle.] (connectivity) for all $R \ge R_0$ and
    $x \in \Omega$, there is a subgraph $\mathcal{\bar{G}}$ in
    Assumption~\ref{ass:network}.2 with
    \begin{equation*}
      \vol(\mathcal{\bar{N}}) \le \nu_1 \left(\frac{R}{R_0}\right)^d,
    \end{equation*}
    and for which the following $d$-dimensional isoperimetric inequality holds,
    \begin{equation*}
      \left(\vol(X)\right)^{(d-1)/d} \le \nu_2 |\mathcal{\bar{E}}(X, \mathcal{\bar{N}} \setminus X)|
    \end{equation*}
    for all $X \subset \mathcal{\bar{N}}$ assuming
    $\vol(X) \le \vol(\mathcal{\bar{N}} \setminus X)$
  \end{enumerate}
  then Lemma~\ref{lem:poincare} holds with $\mu$ depending only on
  $\nu_1$ and $\nu_2$.
\end{lemma}
\begin{proof}
  The eigenvalues $\lambda_1$ and $\lambda_2$ in the proof of
  Lemma~\ref{lem:poincare} for the Dirichlet and Neumann cases,
  respectively, can be related to the $d$-isoperimetric constant for
  the subgraph $\mathcal{\bar{G}}$. For brevity, both cases are
  treated simultaneously. Let $V_1 = \{0\}$ and $V_2 = \spanop(1)$ (containing all constant functions) be
  subspaces of $V(\mathcal{\bar{N}})$ and $\hat V(\mathcal{\bar{N}})$,
  respectively.

We make use of bounds for the normalized graph Laplacian as presented
in \cite{ChungYauSobolev}.  Continuing from \eqref{eq:minmax_dirichlet}
and \eqref{eq:minmax_neumann} above and using that (by assumption)
all edge lengths are smaller than $R_0$, we can obtain the following
bound for $k = 1$ (the Dirichlet case) and $k = 2$ (the Neumann case),
\begin{equation*}
  \begin{aligned}
    \lambda_k & = \min_{\substack{v \ne 0}} \max_{c \in V_k} \frac{(\bar L v, v)}{(\bar M (v - c), v - c)} \\
    & \ge R_0^{-2} \min_{\substack{v \ne 0}} \max_{c \in V_k} \frac{\sum_{x \in \mathcal{\bar{N}}} \sum_{x \sim y} (v(x) - v(y))^2}{\sum_{x \in \mathcal{\bar{N}}} \bar{d}_x (v(x) - c)^2} \\
    & = R_0^{-2} \hat \lambda_k\\
  \end{aligned}
\end{equation*}
where $\hat \lambda_k$ is the $k$th smallest eigenvalue of the
normalized graph Laplacian (see e.g.\ \cite{ChungBook}) in the two cases, and $x \sim y$ means
that $x$ and $y$ are adjacent in $\mathcal{\bar{G}}$.

The following bound
\begin{equation*}
  \hat \lambda_k \ge C_{\nu_2}\vol(\mathcal{\bar{G}})^{-2/d}
\end{equation*}
on the eigenvalues of the normalized graph Laplacian is proved for
$k=1$ (the Dirichlet case) in \cite[Propositon 7.1]{ChungGrigYau00}
and for $k=2$ (the Neumann case) in \cite[Theorem 4]{ChungYauSobolev}.
Using the assumed volume bound, we get
\begin{equation*}
  \lambda_k^{-1/2} \le R_0 \hat \lambda_k^{-1/2} \le C_{\nu_2}^{-1/2} R_0 \vol(\mathcal{\bar{G}})^{1/d} \le C_{\nu_1, \nu_2} R
\end{equation*}
which means we can prove Lemma~\ref{lem:poincare} with
$\mu = R^{-1}\lambda_k^{-1/2} \le C_{\nu_1, \nu_2}$.
\end{proof}
\begin{example}[Grid network]
  For simple graphs, for example a grid, it is possible to bound
  $\nu_1$ and $\nu_2$ analytically and thus provably establish a small
  $\mu$. Let $\mathcal{G}$ be an $n \times n$ grid network in $d=2$
  with nodes in $(i/n, j/n)$ for $i,j=0,\ldots,n$ and edges between
  the nodes which are at distance $1/n$ from each other. Then with
  $R_0 = 1/n$, we can set $\nu_1 = 4$ and $\nu_2 = 4$ (see
  \cite{Tillich}) to satisfy the assumptions of
  Lemma~\ref{lem:isoperimetric}.
\end{example}

\section{A spatial network solver based on subspace decomposition}

We will use a preconditioned conjugate gradient method to solve the model problem. The preconditioner is based on the method proposed in \cite{KY16,KPY18} for elliptic partial differential equations which has its foundation in the rich literature on subspace decomposition and correction, see \cite{Xu90} and references therein. For completeness, we include the full error analysis even though some parts, not directly depending on the underlying network formulation, are classical results.

\subsection{Subspace decomposition}
We make the following subspace decomposition
$$
\VV_0=\VV_H,\quad \VV_j=\VV(U(y_j)),\quad j=1,\dots,m,
$$
and let $\PP_j:\VV\rightarrow \VV_j$ be orthogonal projections fulfilling
\begin{equation}\label{eq:defPj}
(\KK \PP_j \vv,\vv_j)=(\KK \vv,\vv_j),
\end{equation}
for all $\vv_j\in \VV_j$. The existence and uniqueness of such operators $\PP_j$ follows by Assumption \ref{ass:K}.3 and the fact that $\VV_j$ is a subspace of $\VV$. Following \cite{KY16} we introduce the operator 
$$
\PP=\PP_0+\PP_1+\dots+\PP_m.
$$
The operator $\PP$ involves direct solution of decoupled local linear systems and one coarse scale linear system. Therefore, the PCG method we propose below can be referred to as semi-iterative. Under Assumptions \ref{ass:K} and \ref{ass:network} we can prove the following stability of the subspace decomposition, which is essential for the convergence of the method.

\begin{lemma}\label{lem:K1K2}
If Assumptions \ref{ass:K} and \ref{ass:network} hold and $H\geq 2R_0$ then
there is a decomposition $\vv=\sum_{j=0}^m\vv_j$ that satisfies
\begin{equation}\label{eq:C1}
  \sum_{j=0}^m |\vv_j|_{\KK}^2 \leq C_1 |\vv|_{\KK}^2.
\end{equation}
Moreover, every decomposition $\vv=\sum_{j=0}^m \vv_j$ with $\vv_j\in \VV_j$ satisfies
\begin{equation}\label{eq:C2}
  |\vv|_{\KK}^2 \leq C_2 \sum_{j=0}^m |\vv_j|_{\KK}^2.
\end{equation}
The constants are $C_1=C_d\beta\alpha^{-1}\sigma\mu^2$ and $C_2=C_d\beta\alpha^{-1}$, where $C_d$ only depends on $d$.
\end{lemma}
\begin{proof}
Section \ref{sec:coarse_scale_representation} is devoted to the proof of this Lemma.
\end{proof}

Given Lemma \ref{lem:K1K2}
we can bound the operator norm of polynomials of $\PP$ in terms of the spectrum of $\PP$, using the spectral theorem of symmetric operators in finite dimensional spaces. Equations (\ref{eq:C1}--\ref{eq:C2}) give us bounds of the spectrum from above and below $\lambda\in [C_1^{-1},C_2]$. This result is classical and available in other publications, see \cite{KY16} and references therein. We include a proof in the appendix for completeness of the presentation.
\begin{lemma}\label{lem:plambda}
If Assumptions \ref{ass:K} and \ref{ass:network} hold and $H\geq 2R_0$, then the spectrum of $\PP$ fulfills 
$$
C_1^{-1}\leq \lambda\leq C_2
$$ 
and for any polynomial $p$ it holds
\begin{equation*}
  \sup_{\vv \in \VV} \frac{|p(\PP) \vv|_{\KK}}{|\vv|_{\KK}}\leq \max_{\lambda\in[C_1^{-1},C_2]}|p(\lambda)|.
\end{equation*}
\end{lemma}

\subsection{Preconditioned conjugate gradient method}
 The preconditioned conjugate gradient method (PCG) applies the conjugate gradient method to a preconditioned system: find $\uu\in \VV$ such that for all $\vv\in\VV$,
 $$(\BB\KK\uu,\vv)=(\BB\ff,\vv).$$
 In our case $\PP=\BB\KK:\VV\rightarrow \VV$ for some operator $\BB:\hat\VV\rightarrow \VV$ that is not explicitly formed. We have that 
 \begin{equation}\label{eq:TTsym}
(\KK \uu,\PP\vv)=\sum_{j=0}^m(\KK\uu,\PP_j \vv)=\sum_{j=0}^m(\KK\PP_j\uu,\PP_j \vv)=\sum_{j=0}^m(\KK\PP_j\uu,\vv)=(\KK\PP\uu,\vv),
\end{equation}
i.e. $\PP$ is symmetric with respect to the bilinear form induced by $\KK$, and positive definite since $(\KK\PP\vv,\PP\vv)=|\PP\vv|^2_\KK\geq 0$ with equality if and only if  $\vv=0$. 

 From the classical analysis of the PCG algorithm we have that the error in each iteration can be written as
$$
\uu-\uu^{(\ell)}=p_\ell(\PP)(\uu-\uu^{(0)}),
$$
where $p_\ell$ is a polynomial of degree $\ell$ fulfilling $p_\ell(0)=1$, $\uu^{(\ell)}$ is the PCG approximation after $\ell$ steps in the algorithm and $\uu^{(0)}$ is an initial guess. In each iteration the error in the conjugate gradient method $|\uu-\uu^{(\ell)}|_\KK$ is minimized over the Krylov subspace:
$$
\text{span}\{\mathbf{s}^0,\PP\mathbf{s}^0,\PP^2\mathbf{s}^0,\dots,\PP^{\ell-1}\mathbf{s}^0\},
$$
where $\mathbf{s}^0 =\PP( \uu - \uu^{(0)})$. Therefore, the PCG solution realizes the minimum
\begin{equation}\label{eq:boundp}
|\uu-\uu^{(\ell)}|_\KK\leq \min_{\tiny{\begin{array}{c}\text{deg}(p)\leq \ell\\ p(0)=1\end{array}}}|p(\PP)|_\KK|\uu-\uu^{(0)}|_\KK\leq \min_{\tiny{\begin{array}{c}\text{deg}(p)\leq \ell\\ p(0)=1\end{array}}}\max_{\lambda\in[C_1^{-1},C_2]}|p(\lambda)||\uu-\uu^{(0)}|_\KK,
\end{equation}
where we have applied Lemma \ref{lem:plambda} in the last step. 
Finding that minimizing polynomial is a classical min-max problem and the solution is given by a shifted and scaled Chebyshev polynomial, see for instance \cite{CGCheb}. We conclude
\begin{equation}\label{eq:optimal}
|\uu-\uu^{(\ell)}|_\KK\leq 2\left(\frac{\sqrt{\kappa}-1}{\sqrt{\kappa}+1}\right)^\ell|\uu-\uu^{(0)}|_\KK,
\end{equation}
where $\kappa = C_1C_2$.

We are ready to formulate the main convergence result for the proposed preconditioned conjugate gradient method.
\begin{theorem}\label{thm:cg}
If Assumptions \ref{ass:K} and \ref{ass:network} hold and $H\geq 2R_0$ then the preconditioned conjugate gradient approximation $\uu^{(\ell)}$ fulfills
$$
|\uu-\uu^{(\ell)}|_\KK\leq 2\left(\frac{\sqrt{\kappa}-1}{\sqrt{\kappa}+1}\right)^\ell|\uu-\uu^{(0)}|_\KK,
$$
where $\sqrt{\kappa} = \sqrt{C_1C_2}=C_d\beta\alpha^{-1}\sigma^{1/2}\mu$.
\end{theorem}

\section{Stability of the subspace decomposition}
\label{sec:coarse_scale_representation}

This section is dedicated to the proof of the stability of the subspsace decomposition stated in Lemma~\ref{lem:K1K2}.
We first define an interpolant onto the space $V_H$ and then prove an interpolation error bound and a bound for products of functions in $V_H$ and $V$. Together these two results are used to prove Lemma \ref{lem:K1K2}.  

\subsection{Interpolation error bound}

We introduce an interpolation operator from $V$ to
$V_{H}$. For each basis function $\varphi_k$ we denote the
unique element that contains $y_k$ by $T_k$ and let 
$\psi_k=|1|^{-2}_{M,T_k}$ be a constant function.
The interpolation operator is
then defined by
$$\II v = \sum_{k=1}^{m_0} (M_{T_k} \psi_k, v) \varphi_k.$$
For future reference we note that we have the bound 
\begin{equation}\label{eq:interpolation_stability}
|\psi_k|_{M,T_k}=|1|_{M,T_k}^{-1}\leq \rho^{-1/2}H^{-d/2}
\end{equation}
for all $H\geq 2R_0$ by Assumption \ref{ass:network}.
  The idea of the construction of $\mathcal{I}$ comes from
  the construction of the 
  Cl\'{e}ment finite element interpolant \cite{C75}.
  There is some freedom in the choice of domain
  for 
  $\psi_k$.  We picked $T_k$ for simplicity. Next, we present some
  auxiliary lemmas needed in the proof of the interpolation bound.

\begin{lemma}
  \label{lem:bound_norm}
  If Assumption~\ref{ass:network} holds and $H\geq 2R_0$, then 
  \begin{equation}
    |\varphi_k|_{M,T} \le \sigma^{1/2}\rho^{1/2}H^{d/2}
  \end{equation}
  for mesh nodes $k = 1,\ldots,m$. If $T \cap U(y_k)$ is empty, we have $|\varphi_k|_{M,T} = 0$.
\end{lemma}
\begin{proof}
  Directly from Assumption~\ref{ass:network}.1 we have 
  \begin{equation*}
    |\varphi_k|_{M,T}^2 \le |1|_{M,T}^2 \le \sigma \rho H^{d}.
  \end{equation*}
  Further, since the support of $\varphi_k$ is a subset of $U(y_k)$ (the elements adjacent to node $y_k$)
  and $|v|_{M,T}$ depends only on values of $v$ in $T$, we conclude
  that $|\varphi_k|_{M,T} = 0$ if $T \cap U(y_k)$ is empty.
\end{proof}

\begin{lemma}
  \label{lem:bound_energy_norm}
  If Assumption~\ref{ass:network} holds and $H\geq 2R_0$, then 
  \begin{equation}
     | \varphi_k |_{L,T} \le \sigma^{1/2} \rho^{1/2} H^{d/2-1}
  \end{equation}
  for mesh nodes $k = 1,\ldots,m$ and $H\geq 2R_0$. If $U(T) \cap U(y_k)$ is empty, we have $| \varphi_k |_{L,T} = 0$.
\end{lemma}
\begin{proof}
  We consider an arbitrary $T \in \mathcal{T}_H$ and observe that the
  basis functions $\varphi_k$ have Lipschitz constant $H^{-1}$. Using
  this, and the network properties from Assumption~\ref{ass:network}.1,
  we get
  \begin{equation}
    \label{eq:bound_energy_norm}
    \begin{aligned}
      | \varphi_k|_{L,T}^2 & = \sum_{x \in \mathcal{N}(T)} (L_x \varphi_k, \varphi_k)
      = \frac{1}{2} \sum_{x \in \mathcal{N}(T)} \sum_{x \sim y}  \frac{(\varphi_k(x) - \varphi_k(y))^2}{|x - y|} \\
      & \le \frac{H^{-2}}{2} \sum_{x \in \mathcal{N}(T)} \sum_{x \sim y} |x - y|
      = H^{-2} |1|_{M,T}^2 \le \sigma \rho H^{d-2}.
    \end{aligned}
  \end{equation}
  To see that the norm is zero when $U(T) \cap U(y_k)$ is empty, we
  note that $|\varphi_k|_{L,T}$ depends on values of $\varphi_k$ at
  nodes adjacent to nodes in $T$, which by
  Assumption~\ref{ass:network}.3 is a subset of the nodes in
  $U(T)$. Since the support of $\varphi_k$ is $U(y_k)$ and the
  intersection between $U(T)$ and $U(y_k)$ is empty, the norm is zero.
\end{proof}

With the above results and assumptions, we establish an interpolation bound. The constants $C_d$ in the calculations below only depends on the dimension $d$ and can change between equations. 
\begin{lemma} 
  \label{lem:interpolation}
  If Assumption~\ref{ass:network} holds and $H\geq 2R_0$, then for $v \in V$,
  \begin{equation}
   H^{-1}| v - \II v |_{M} + |\II v|_{L} \le C_{d} \sigma^{1/2}\mu|v|_{L}.
  \end{equation}  
\end{lemma}
\begin{proof}
  Consider a single element $T\in \mathcal{T}_H$. Without loss of
  generality, we order the free mesh nodes (i.e.\ the $m_0$ first
  nodes) so that nodes with indices $k = 1, \ldots, \hat m_0$ are
  those for which $U(T) \cap U(y_k) \ne \emptyset$. By our choice of a
  regular mesh we note that $\hat m_0 \le 3^d$. With this ordering, we
  have from Lemma~\ref{lem:bound_energy_norm} that
  $|\varphi_k|_{L,T} = 0$ whenever $k > \hat m_0$. Using this and the
  definition of $\II$ (which can also be applied to functions in
  $\hat V$), we have for any constant function $c \in \hat V$,
  \begin{equation}
    \begin{aligned}
      \label{eq:I_definition}
      |\II (v - c)|_{L,T} & \le \left| \sum_{k=1}^{\hat m_{0}}(M_{T_k} \psi_k, v - c) \varphi_k \right|_{L,T} \\
      & \le \sum_{k=1}^{\hat m_0} | \psi_k |_{M,T_k} | v - c |_{M,T_k} | \varphi_k |_{L,T} \\
      & \le | v - c |_{M,U_2(T)} \sum_{k=1}^{\hat m_0} | \psi_k |_{M,T_k} |\varphi_k|_{L,T} \\
      & \le C_d\sigma^{1/2} H^{-1} | v - c |_{M,U_2(T)},
    \end{aligned}
  \end{equation}
  where we used equations~\eqref{eq:interpolation_stability} and
  \eqref{lem:bound_energy_norm} in the last step. We distinguish
  between two cases.

  If $U(T)$ does not intersect with $\Gamma$, then
  $|c - \II c|_{L,T} = 0$, and by the Poincaré inequality in
  Lemma~\ref{lem:poincare} (with $R = 5H/2$ and $x$ such that
  $B_R(x) = U_2(T)$) we can choose $c$ and bound
  $| v - c |_{M,U_2(T)} \le \mu H |v|_{L,U_3(T)}$. We conclude
  \begin{equation}\label{eq:U3}
    \begin{aligned}
      | v - \II v |_{L,T} & \le | v - c |_{L,T} + |\II(v - c)|_{L, T} \\
       & \le C_d\sigma^{1/2}\mu |v|_{L, U_3(T)}. \\
    \end{aligned}
  \end{equation}
  
  Otherwise, if $U(T)$ does intersect with $\Gamma$, then so does $U_2(T)$.
  Since $U_2(T)$ is a patch of disjoint elements that intersect with
  the boundary $\Gamma$, and since we have assumed that $\Gamma$ is a
  union of element faces, $U_2(T)$ contains an element face contained
  in $\Gamma$. Since the faces are boxes in dimension $d-1$ of radius
  $R \ge R_0$, then by Assumption~\ref{ass:network}.4, the face, and
  consequently $U_2(T)$, contains at least one boundary node. Now,
  \eqref{eq:U3} holds again, but with $c = 0$ using the Friedrichs
  inequality in Lemma~\ref{lem:poincare} on the box $B_R(x) = U_2(T)$.

  By taking squares in \eqref{eq:U3} and summing over all elements we
  conclude $ | v - \II v |_{L} \le C_d\sigma^{1/2}\mu |v|_{L}$.  With
  the same argument, but starting with the $|\cdot|_{M,T}$ norm in
  \eqref{eq:I_definition} and using Lemma~\ref{lem:bound_norm} instead
  of Lemma~\ref{lem:bound_energy_norm} we have
  $| v - \II v |_{M} \le C_d\sigma^{1/2}\mu H |v|_{L}$. Finally,
  $$| \II v |_{L} \le | v - \II v |_{L} + |v|_{L} \le C_{d}\sigma^{1/2}\mu| v |_{L}$$ concludes the proof.
\end{proof}

We also need to bound the product of a basis function and an arbitrary function.
\begin{lemma}
  \label{lem:product_rule}
  It holds that
  \begin{equation}
     | v \varphi_k |_{L,T}^2 \le 2 \left(H^{-2} | v |_{M,T}^2 + | v |_{L,T}^2\right)
  \end{equation}
  for all mesh nodes $k = 1,\ldots,m$ and $v \in \hat V$.
\end{lemma}
\begin{proof}
  We drop the subscript $k$ and call $\varphi = \varphi_k$. By using
  that $\varphi$ takes values only between 0 and 1, and that the
  Lipschitz constant of $\varphi$ is $H^{-1}$, we get
  \begin{equation*}
    \begin{aligned}
      |v \varphi|_{L,T}^2 &= (L_T (v \varphi), v \varphi) = \frac{1}{2} \sum_{x \in \mathcal{N}(T)} \sum_{x \sim y} \frac{\left(v(x)\varphi(x) - v(y)\varphi(y)\right)^2}{|x - y|} \\
      & = \frac{1}{2} \sum_{x \in \mathcal{N}(T)} \sum_{x \sim y} \frac{\left(v(x)(\varphi(x) - \varphi(y)) + (v(x) - v(y))\varphi(y)\right)^2}{|x - y|} \\
      & \le \sum_{x \in \mathcal{N}(T)} \sum_{x \sim y} \frac{v(x)^2 |x - y|^2 H^{-2}+ (v(x) - v(y))^2}{|x - y|} \\
      & = 2\left(H^{-2}|v|_{M,T}^2 + |v|_{L,T}^2\right).\\
    \end{aligned}
  \end{equation*}
\end{proof}

\subsection{Proof of Lemma \ref{lem:K1K2}}
Under Assumptions \ref{ass:K} and \ref{ass:network} we want to prove that 
there is a particular decomposition $\vv=\sum_{j=0}^m \vv_j$ that satisfies
\begin{equation}\label{eq:trirev}
  \sum_{j=0}^m |\vv_j|_{\KK}^2 \leq C_d\beta\alpha^{-1}\sigma\mu^2 |\vv|_{\KK}^2,
\end{equation}
and that every decomposition $\vv=\sum_{j=0}^m \vv_j$ satisfies
\begin{equation}\label{eq:tri}
  |\vv|_{\KK}^2 \leq C_d\beta\alpha^{-1} \sum_{j=0}^m |\vv_j|_{\KK}^2.
\end{equation}
To prove the first bound we use the decomposition $$\vv_0=\III\vv\quad\text{and}\quad\vv_j=(\vv-\III \vv)\varphi_j\in \VV_j$$
for $j=1,\dots,m$.

\begin{proof}[Proof of Lemma \ref{lem:K1K2}]
We start with equation (\ref{eq:tri}) in the $\LL$-norm. By construction
$$
 |\vv|^2_{\LL}\leq 2|\vv_0|^2_\LL+2\left|\sum_{j=1}^m\vv_j\right|^2_\LL.
$$
To bound the second term we pick a $T \in \mathcal{T}_H$. Since
  $\vv_j \in \VV(U(y_j))$ and $\LL_T \vv = 0$ for
  $\vv \in \VV(\Omega \setminus U_2(y_j))$, we have that $\LL_T \vv_j$ can
  be non-zero for at most $C_d$ mesh nodes $j$, where $C_d$ depends
  only on $d$. Since $\LL_T$ is local in this sense, we get
  \begin{equation*}
    \left|\sum_{j=1}^m\vv_j\right|^2_{\LL,T} \le C_d \sum_{j=1}^m |\vv_j|_{\LL,T}^2.
  \end{equation*}
  Summing over $T \in \mathcal{T}_H$ proves the inequality in
  $\LL$-norm
$$
 |\vv|_{\LL}^2\leq 2|\vv_0|_\LL+2C_d \sum_{j=1}^m |\vv_j|_{\LL}^2\leq 2C_d \sum_{j=0}^m |\vv_j|_{\LL}^2.
$$  
  From Assumption~\ref{ass:K} the global $\LL$- and
  $\KK$-norms are equivalent with equivalence constant
  $\beta\alpha^{-1}$. We thus get the asserted inequality with
  $C_2 = C_{d}\beta\alpha^{-1}$.

  To prove equation (\ref{eq:trirev}), we use Lemma~\ref{lem:product_rule} and \ref{lem:interpolation}
  componentwise, and a similar locality
  argument of $\LL_T$ as in the previous paragraph and get
  \begin{equation*}
    \begin{aligned}
      \sum_{j=1}^m |\vv_j|_{\LL}^2 & = \sum_{j=1}^m |(\vv-\III\vv)\varphi_j|_{\LL}^2 \\
      & \le 2 \sum_{j=1}^m \sum_{T \subset U_2(y_j)} \left(H^{-2}|\vv-\III\vv|_{\MM,T}^2 + |\vv-\III\vv|_{\LL,T}^2\right) \\
      & \le C_d  \left(H^{-2}|\vv-\III\vv|_{\MM}^2 + |\vv-\III\vv|_{\LL}^2\right) \\
      & \le C_{d}\sigma\mu^2 |\vv|_{\LL}^2.
    \end{aligned}
  \end{equation*}
  Furthermore, $|\vv_0|_\LL= |\III\vv|_\LL\leq C_{d}\sigma^{1/2}\mu|\vv|_\LL$. Altogether we have 
  $$
  \sum_{j=0}^m |\vv_j|_{\LL}^2\leq C_{d}\sigma\mu^2 |\vv|_{\LL}^2.
  $$
    Again, we use equivalence of $\LL$- and $\KK$-norms to get the
  first inequality \eqref{eq:trirev} with $C_d\beta\alpha^{-1}\sigma\mu^2$.
\end{proof}

\section{Numerical examples}\label{sec:numexp}
We start by investigating the network assumptions with particular focus on homogeneity and connectivity for a number of spatial network models. We then solve heat conductivity and structural problems using the proposed semi-iterative scheme and study the convergence rate.

\subsection{Network assumptions} \label{sec:network_analysis}

\begin{figure}[b]
  \centering
  \includegraphics{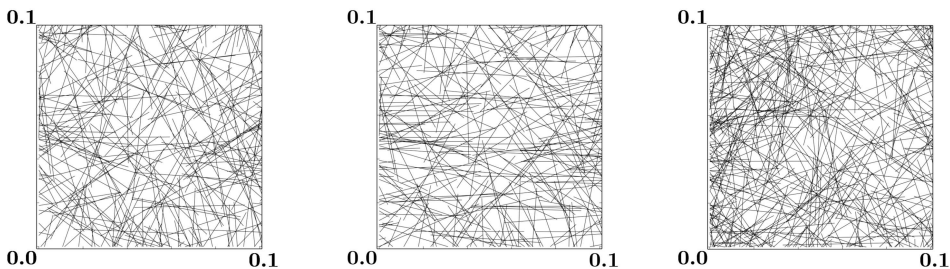}
  \caption{The lower left corners ($[0.0,0.1]^2$) of the three types of fiber networks analyzed. The left network was generated completely uniformly, the center introduce a bias in fiber rotation, and the right a bias in fiber placement.}
  \label{fig:networks}
\end{figure}

\begin{figure}
  \centering
  \includegraphics{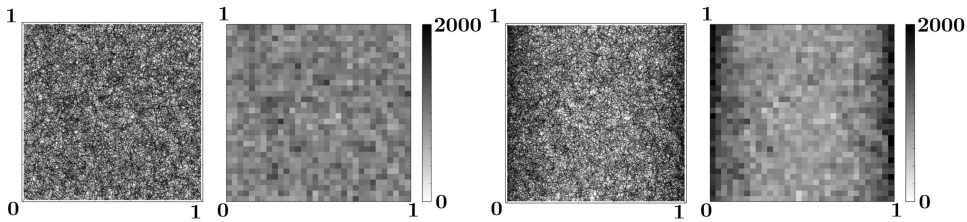}
  \caption{Entire networks and illustrations of the homogeneity assumption for the uniform fiber network (left) and the network with more fibers along the lines $x_1=0$ and $x_1=1$ (right). The boxes represent $(2R)^{-d}|1|^2_{M,B_R(x)}$ for the respective box $B_R(x)$. }
  \label{fig:randHomo}
\end{figure}

In the following example, the connectivity and homogeneity assumptions in Assumption~\ref{ass:network} are visualized by generating and analyzing three types of random two-dimensional fiber networks. These networks represent different forms of irregularity in the network structure. The first is created by placing and rotating fibers uniformly in the domain, the second network has a bias in the fiber rotation, and the third has a bias in the fiber placement. The lower left corners of these three networks are presented in Figure~\ref{fig:networks}, and illustrations of entire networks can be found in Figure~\ref{fig:randHomo}.

All three networks are generated by placing fibers with a fixed length ($r=0.05$) in the unit square until a specified density is achieved. Density in this example is defined as the total edge length, $|1|_M^2 = 1000$. Each fiber is placed randomly with its midpoint in $[-0.5r,1+0.5r]^2$, then rotated given a random degree. Any part outside the unit square is removed. The slightly extended domain is for uniformity along the boundaries, and in this example, networks are generated with different distributions for the midpoint placement and angle of the fibers. Then every two intersecting fibers are connected by placing a node in each intersection. After this step, any two nodes closer than $r\times 10^{-4}$ are merged, setting a lower bound on the edge lengths. The final step is to remove any potential disconnected edges that have formed along the boundary. Networks generated with these parameters have around 0.3 million network nodes. 

To illustrate the homogeneity assumption, $(2R)^{-d}|1|^2_{M,B_R(x)}$ is presented for the networks discussed. For these networks, five different regular grids are considered; $R^{-1} = 4, \ 8, \ 16, \ 32, $ and $64$. For each element in the grids and network types, the corresponding expression $(2R)^{-d}|1|^2_{M,B_R(x)}$ is evaluated. The evaluations for the uniformly distributed network and the network with biased midpoint placement are illustrated in Figure~\ref{fig:randHomo}, with the numerical results for all the setups presented in Table~\ref{table:eig_dat}

The connectivity assumption in Assumption~\ref{ass:network} is analyzed using subgraphs. The assumption holds if we can find a connected subgraph, $\overline{\mathcal{G}}$, on $B_{R+R_0}(x)$. This subgraph has to contain all nodes and edges containing nodes in $B_R(x)$, where the second-smallest eigenvalue of the corresponding eigenvalue problem, $\overline{L} u = \lambda_2 \overline{M} u,$ has the property: $\lambda_2^{-1} \leq \mu^2R^2$, given some constant $\mu>0$ (as discussed in the proof of Lemma~\ref{lem:poincare}). For simplicity, we only find one such network for each inequality considered, where the network found has a minimal amount of edges using a breadth-first search scheme. 

\begin{figure}
\centering
  \includegraphics{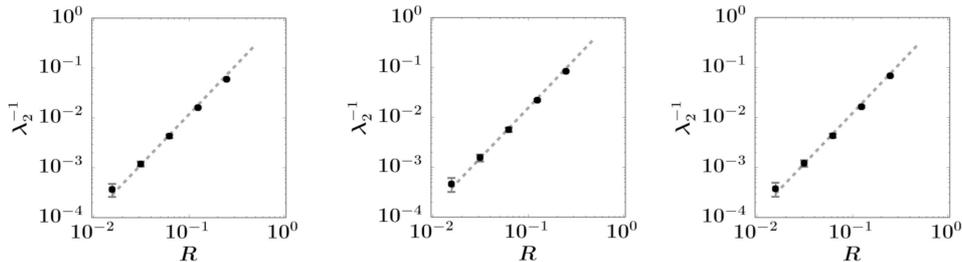}
  \caption{The relation between $\lambda_2^{-1}$ and $R$ for the following networks: uniform (left),  fiber orientation bias (center), and fiber placement bias (right). The dots represent the average $\lambda_2^{-1}$ for the given grid size, and the feet the standard deviation.}
  \label{fig:randConn}
\end{figure}

Similar to the homogeneity example, the unit square is partitioned into a regular grid of squares. The inequality induced by each square, $B_R(x)$, with length scale $R_0^{-1}=64$ for each network is evaluated. In Figure~\ref{fig:randConn}, the $\lambda_2^{-1} \leq \mu^2R^2$  inequality (connected to the Poincar\'{e} inequality) is illustrated. The upper bounds on $\mu$ observed from each grid considered in the connectivity experiments are presented in Table~\ref{table:eig_dat}. 
\begin{table}[h]
\centering
\caption{Table with $(\sigma,\mu)$ for different $R$ attained through numerical evaluation.}
\begin{tabular}{l|c|c|c|c|c|c}
\ & $R^{-1} = 4$ & $R^{-1} = 8$  & $R^{-1} = 16$ & $R^{-1} = 32$ & $R^{-1} = 64$ \\
\hline 
Uniform   & (1.04,0.49)& (1.08,0.53) & (1.27,0.57) & (1.85,0.675) & (3.42,1.53)  \\
Rand. Orient.  & (1.04,0.59)  & (1.08,0.61) & (1.27,0.69) & (1.87,0.83) & (2.93,1.35)   \\
Rand. Domain  & (1.04,0.53) & (1.57,0.54)  & (2.13,0.58)& (3.1,0.76) & (6.86,1.45)  
\end{tabular}
\label{table:eig_dat}
\end{table}

The results illustrate that $\lambda_2^{-1}\sim R^{-2}$ for all considered discretization levels, meaning that the considered networks fit the connectivity assumptions down to $R^{-1}=64$. For the worst case scenario considered in Table~\ref{table:eig_dat} we detect an increase in $\sigma$ and $\mu$ with decreasing $R$ as expected. It is most evident in the finest grid size. It is observable that the network with bias in the fibers' orientation has lower connectivity than the uniformly distributed network but has similar homogeneity. The network with bias in the fiber midpoint placement has lower homogeneity but comparable connectivity to the non-biased network. Similar results also hold for the corresponding eigenvalue $\lambda_1$ related to the Friedrich inequality.  

\subsection{Heat conductivity}
\begin{figure}
  \centering
  \includegraphics{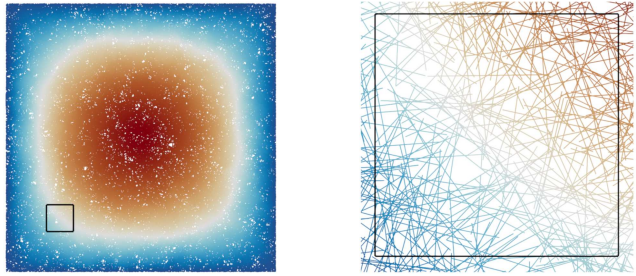}
  \caption{Illustrations of the solution to a weighted graph Laplacian problems on the form \eqref{eq:heat}. The left picture shows the entire solution, and the right shows the highlighted area in detail. }
  \label{fig:lap_sol}
\end{figure}

We consider a heat conductivity problem for three different two-dimensional networks ($d=2$). The solutions can be interpreted as the (scalar) temperature in each node, with the operator $K$ defined in equation \eqref{eq:scalarmodel}. We seek $u\in V=\{v\in \hat V\,:\,v(x)=0,x\in \partial\Omega \}$ such that for all $v\in V$,
\begin{equation} \label{eq:heat}
  (K u,v) = (M1,v).
\end{equation}
We let  $f=M1$ in order to distribute the source equally over the edges of the network. An illustration of a solution to such a problem is presented in Figure~\ref{fig:lap_sol}.

\begin{figure}
  \centering
  \includegraphics{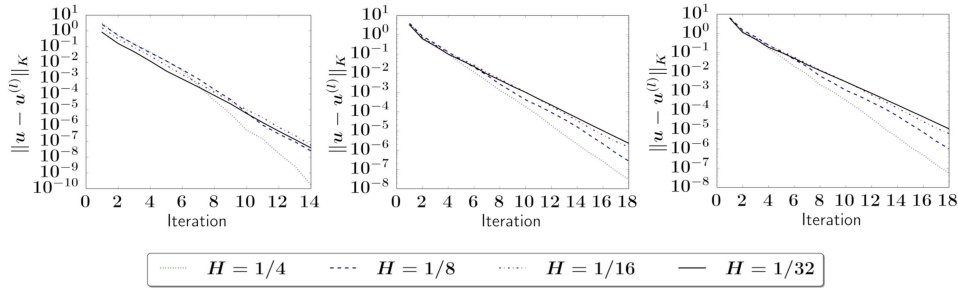}
  \caption{The error for each iteration of the three weighted graph Laplacians analyzed. From left to right: constant weights and a uniform regular grid as the network, a fiber based network with constant weights, and a fiber based network with random weights.}
  \label{fig:lap_conv}
\end{figure}

The networks considered are a uniform regular grid with $(2^{9}+1)^2$ nodes, a uniform fiber based network as introduced in Section~\ref{sec:network_analysis} with constant weight $\gamma_{xy}=1$, and a uniform fiber based network with weights $\gamma_{xy}$ drawn from a uniform distribution on $[0.1,1]$. 
The errors, in $K$-norm, for each iteration with the PCG method is presented in Figure~\ref{fig:lap_conv}. The convergence rates, $\tau_{(\ell)}$:   
\begin{equation} \label{eq:conv_rate}
\tau_{(\ell)}\vcentcolon=\frac{|\uu-\uu^{(\ell-1)}|_\KK}{|\uu-\uu^{(\ell)}|_\KK} \leq \tau, \ l = 2, \ \hdots
\end{equation}
on average, $\overline{\tau}$, and upper bounds, $\tau$, for the three cases can be found in Table~\ref{table:lap_conv}. In the figures and results it is clear that the method converges exponentially for all three cases, which is consistent with Theorem~\ref{thm:cg}.

\begin{table}[h]
\centering
\caption{The PCG method's average and worst convergence rates, $(\overline{\tau}, \tau)$, for different heat conductivity problems.}
\begin{tabular}{l|c|c|c|c}
Problem & $H^{-1} = 4$  & $H^{-1} = 8$ & $H^{-1} = 16$ & $H^{-1} = 32$ \\
\hline 
Grid            & (0.18,0.31) & (0.25,0.33) & (0.27,0.32) &	(0.28,0.31)   \\
Fiber           & (0.29,0.48) &	(0.33,0.43) & (0.39,0.49) &	(0.42,0.47)   \\
Fiber weighted  & (0.34,0.44) & (0.40,0.52) & (0.45,0.51) & (0.47,0.54) 
\end{tabular}
\label{table:lap_conv}
\end{table}

\subsubsection{Analysis of the constant $C_d$}
In the following discussion, the theoretical bound on the convergence rate defined by Theorem~\ref{thm:cg} is evaluated closer. Given the theorem and \eqref{eq:conv_rate} we have the following bound:
\begin{equation} \label{eq:conv_rate_2}
\tau = \frac{\sqrt{\kappa}-1}{\sqrt{\kappa}+1}, \ \sqrt{\kappa} =  C_d\alpha^{-1}\beta\sqrt{\sigma}\mu.
\end{equation}

With the same approach to the discussion in Section \ref{sec:network_analysis}, values of $\sigma$ and $\mu$ for the uniform regular grid network are numerically evaluated. Together with the convergence results (average) in Table \ref{table:lap_conv}, and the fact that $\alpha = \beta = 1$ for this example, the constant $C_d$ in Theorem~\ref{thm:cg} can be approximated by
$$C_d = \frac{1+\overline\tau}{\sqrt{\sigma}\mu(1-\overline\tau)}.$$
The geometric constants and the approximations of $C_d$ are presented in Table~\ref{table:geo_grid}. From these results, it is clear that the constant $C_d\approx 3.5$ bounds the average convergence rates for all the grids analyzed for the regular grid network.

\begin{table}[h]
\centering
\caption{Investigation of the constant $C_d$ using the average convergence rates, $\overline\tau$, and the geometrical constants $\sigma$ and $\mu$ for the regular grid network.}
\begin{tabular}{l|c|c|c|c}
\ & $H^{-1} = 4$  & $H^{-1} = 8$ & $H^{-1} = 16$ & $H^{-1} = 32$ \\
\hline 
$\sigma$  & 1.0	& 1.0 &	1.1 & 1.1    \\
$\mu$  & 0.45 &	0.47 &	0.48 &	0.51 \\
\hline
$C_d$  & 3.2 &	3.5 &	3.5 &	3.2 
\end{tabular}
\label{table:geo_grid}
\end{table}

The convergence results for the grid type network is now compared to the convergence results for the uniformly distributed fiber based network without weights. In Section \ref{sec:network_analysis}, the geometrical constants $\sigma$ and $\mu$ were evaluated for the fiber based network. Using \eqref{eq:conv_rate_2} with the constants in Table~\ref{table:eig_dat}, $C_d=3.5$, and the convergence rates in Table~\ref{table:lap_conv}, the convergence rate of the PCG method on the fiber model is compared to
$$\tilde{\tau} = \frac{\sqrt{\kappa}-1}{\sqrt{\kappa}+1}, \  \sqrt{\kappa} =  3.5\sqrt{\sigma}\mu $$
in Table \ref{table:fib_lap}.

\begin{table}[h]
\centering
\caption{Comparing convergence rate estimates, $\tilde{\tau}$, with average convergence rates, $\overline\tau$, in numerical experiments for a fiber based network.}
\begin{tabular}{l|c|c|c|c}
\ & $H^{-1} = 4$  & $H^{-1} = 8$ & $H^{-1} = 16$ & $H^{-1} = 32$ \\
\hline 
$\overline\tau$   & 0.29 &	0.33 &	0.39 &	0.42
 \\
\hline
$\tilde{\tau}$  & 0.27 & 0.35 & 0.45 & 0.60
\end{tabular}
\label{table:fib_lap}
\end{table}

The convergence estimates in Table \ref{table:fib_lap} are similar to the numerical results except for the last data point. In that case, the convergence rate of the PCG method converges faster than the estimate. With this in mind, it is clear that the convergence rates are affected by the homogeneity and connectivity constants $\sigma$ and $\mu$. Moreover, a similar geometric constant bounds the average convergence rates of the PCG method for the girds analyzed. We conclude that the error is reduced by a constant factor in each iteration, in agreement with Theorem~\ref{thm:cg}. In addition, the bound appears to be sharp in the dependency of $\sigma^{1/2}$ and $\mu$ with a proportionality constant $C_d\approx 3.5$ on coarser scales $H$. 

\subsection{Structural problem for fiber based material}
Here two variations of Example \ref{ex:fiber} is presented, where the networks should be interpreted as an anisotropic mesh of round steel wires with radius $r_w = 2.5\times 10^{-3}$ m. Equation \eqref{eq:tensile} becomes a linearized version of Hooke's law with parameter $\gamma_{xy} =\gamma_1 = AE$, where $A = \pi r_w^2$ is the cross-section area of the wire and $E = 210$\,GPa the wires Young's modulus. The bending forces are handled by the addition of the equations in \eqref{eq:bend}. These additions are linearized versions of Euler--Bernoulli beam theory with parameters $\gamma_{xyz} = EI(|x-y|+|x-z|)^{-2}$ where $E$ is the same Young's modulus and $I = 0.25\pi r_w^4 = 0.25Ar_w^2$ is the second moment of area of the wire. The two coefficients have the following relation,
$$\gamma_{xyz} = EA\frac{r_w^2}{4(|x-y|+|x-z|)^{2}} = \gamma_{1}\frac{r_w^2}{4(|x-y|+|x-z|)^{2}},$$
where $0.05 \leq \frac{r_w}{|x-y|}\leq 500$ for any edge $x\sim y$.

\subsubsection{Pure displacement problem}
\begin{figure}[h]
  \centering
  \includegraphics{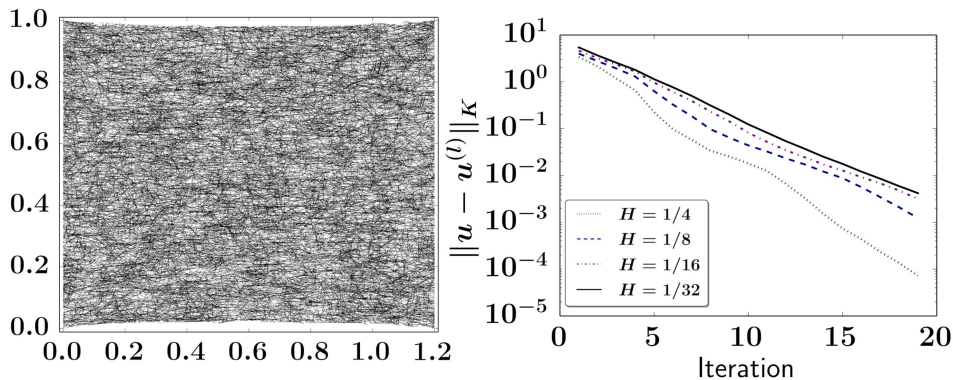}
  \caption{The displaced network along with the convergence results for the PCG method on the two-dimensional structural problem.}
  \label{fig:struct2d}
\end{figure}

The first structural problem is a tensile simulation, where one side is fixed, and the other is displaced, straining the network. For this experiment, the second fiber-based network in Section \ref{sec:network_analysis} is used. This network places the fibers uniformly in the domain but with rotational bias. We only consider forces and displacements in the plane the network resides in for this simulation, meaning that the lateral component is left out. We apply non-homogeneous Dirichlet boundary conditions on the two fixed sides of $\Omega$ and seek $\hat{\uu} = \uu + \gg$ where $\uu \in \VV$ and $\gg(x) = [0.2x_1,0]$. We define $\Gamma = \{x = (x_1, x_2) \in \partial \Omega\,:\, x_1 = 0 \text{ or } x_1 = 1\}$ and solve for $\uu \in \VV = \{\vv \in \hat \VV\,:\, \vv(x)=0, x \in \Gamma \}$ such that for all $v\in V$,
\begin{equation*}
  (\KK \uu, \vv) = -(\KK \gg, \vv).
\end{equation*}
The problem is solved using a direct solver and the PCG method. An illustration of the solution and the errors in each iteration of the PCG method are presented in Figure~\ref{fig:struct2d}, with average convergence rates presented in Table~\ref{tab:struct2D}. The PCG method converges exponentially, and with $\alpha = 0.05$ and $\beta = 500$  the method converges substantially faster than the theoretical upper bound.

\begin{table}[h]
  \centering
  \caption{Average convergence rates, $\overline\tau$, for the two-dimensional structural problem.} \label{tab:struct2D}
  \begin{tabular}{l|c|c|c|c}
    \ & $H^{-1} = 4$  & $H^{-1} = 8$ & $H^{-1} = 16$ & $H^{-1} = 32$ \\
    \hline 
      $\overline\tau$  & 0.54	& 0.62 &	0.65 &	0.65
    \end{tabular}
\end{table}
\subsubsection{Displacement problem with lateral load}
\begin{figure}
  \centering
  \includegraphics{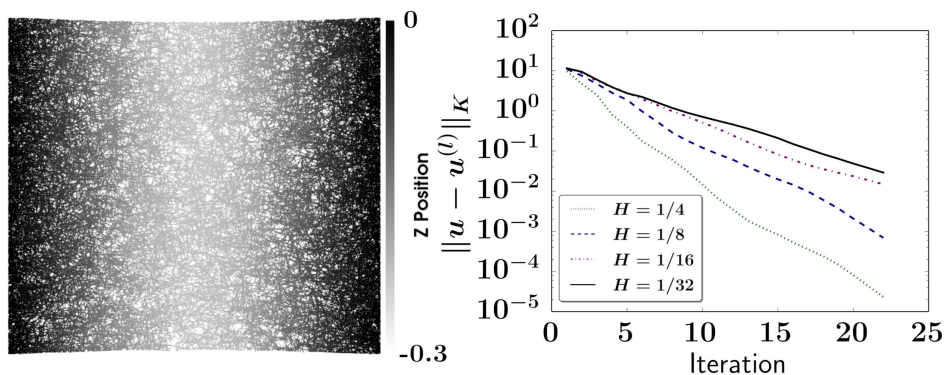}
  \caption{The solution to the exact problem $\mathbf{u}$, along with the convergence results for the PCG method for different $H$.}
  \label{fig:struct3d}
\end{figure}

\begin{table}
  \centering
  \caption{Average convergence rates, $\overline\tau$, for the three-dimensional structural problem.} \label{tab:struct3d}
  \begin{tabular}{l|c|c|c|c}
    \ & $H^{-1} = 4$  & $H^{-1} = 8$ & $H^{-1} = 16$ & $H^{-1} = 32$ \\
    \hline 
    $\overline\tau$ & 0.56 &	0.64 & 0.74 & 0.76
    \end{tabular}
\end{table}
In the second numerical example of the structural problem, we introduce a lateral load to the previous tensile problem with half the strain. This results in a model with $d=2$ and $n = 3$. Using the same $\Gamma$ as in the previous numerical example and $\VV$ defined analogously, the problem is to find $\uu \in \VV$ such that for all $\vv \in \VV$,
\begin{equation*}
  (\KK \uu, \vv) = (\MM \gg' - \KK \gg, \vv) \\
\end{equation*}
where $\gg' = [0,0,-10^3]$ and $\gg(x) = [0.1x_1,0,0]$. The network is the third fiber based network in Section~\ref{sec:network_analysis}. In this network, the fibers are placed with a bias in the domain, resulting in higher densities along the $x_2$-axis. The solution and convergence results for the domain decomposition method for this setup is presented in Figure~\ref{fig:struct3d} and Table~\ref{tab:struct3d}. Similar to the two-dimensional case, the convergence is exponential, and the rates observed are better than the theoretical bound.

\subsubsection*{Acknowledgements.} The first author is supported by the Swedish Foundation for Strategic Research (SSF). The second and third author are supported by the G\"{o}ran Gustafsson Foundation for Research in
Natural Sciences and Medicine and the Swedish Research Council (project number 2019-03517).

\bibliographystyle{amsplain}


\section*{Appendix}
For completeness we include the proof of Lemma \ref{lem:plambda}
\begin{proof}
First we prove the relation
\begin{equation}\label{eq:C1C2}
  C_1^{-1}|\vv|_{\KK}^2 \leq (\KK \PP \vv, \vv) \leq C_2 |\vv|_{\KK}^2
\end{equation}
for all $\vv\in\VV$.

The Cauchy-Schwarz inequality gives
$$
|\vv|_{\KK}^2=\sum_{j=0}^m (\KK \vv,\vv_j)=\sum_{j=0}^m (\KK \PP_j \vv, \vv_j)\leq \left(\sum_{j=0}^m|\PP_j \vv|^2_\KK\right)^{1/2} \left(\sum_{j=0}^m |\vv_j|_\KK^2\right)^{1/2}.
$$ 
Since $\PP_j$ is $\KK$-orthogonal $\sum_{i=0}^m |\PP_j\vv|_\KK^2=\sum_{j=0}^m(\KK \PP_j \vv,\vv)=(\KK \PP\vv,\vv)$
and by Lemma~\ref{lem:K1K2} $\sum_{i=0}^m |\vv_j|_\KK^2\leq C_1|\vv|_{\KK}^2$. We conclude
$
|\vv|_{\KK}^2\leq C_1(\KK \PP\vv,\vv).
$
For the second inequality we again use Lemma~\ref{lem:K1K2} and get
$$
(\KK\PP\vv,\vv)\leq \left|\sum_{j=0}^m\PP_j\vv\right|_\KK |\vv|_\KK\leq \left(C_2\sum_{j=0}^m|\PP_j\vv|^2_\KK\right)^{1/2}|\vv|_\KK\leq \left( C_2 (\KK\PP\vv,\vv) \right)^{1/2}|\vv|_\KK.
$$
Equation (\ref{eq:C1C2}) follows by division by $(\KK\PP\vv,\vv)^{1/2}$ and squaring the equation. 

The spectral theorem for finite dimensional real spaces guarantees the existence of eigenfunctions $\{\zz_i\}_{i=1}^l$ to the $\KK$-symmetric operator $\PP$:
$$
(\KK \PP\zz_i,\vv)=\lambda_i(\KK\zz_i,\vv)
$$
for all $\vv\in \VV$. The eigenfunctions form an orthonormal basis of $\VV$ i.e. $(\KK\zz_i,\zz_j)=\delta_{ij}$.  Furthermore,
using equation (\ref{eq:C1C2}) we get
\begin{equation*}
C_1^{-1} |\zz_i|_\KK^2\le (\KK \PP \zz_i, \zz_i)\leq C_2 |\zz_i|^2_\KK,
\end{equation*}
and conclude
\begin{equation}\label{eq:boundlambda}
C_1^{-1}\leq \lambda_i\leq C_2, 
\end{equation}
for all $i=1,\dots,l$ which proves the first statement of the lemma.
With the eigenfunction expansion $\vv = \sum_{i = 1}^l \alpha_i \zz_i$ the second part follows directly by
\begin{equation*}
\frac{|p(\PP) \vv|^2_{\KK}}{|\vv|_\KK^2}=\frac{\sum_{i,j=1}^l\alpha_i\alpha_j(\KK p(\PP) \zz_i,p(\PP)\zz_j)}{\sum_{i=1}^l \alpha_i^2} 
\end{equation*}
\begin{equation*}
=\frac{\sum_{i=1}^l\alpha_i^2p(\lambda_i)^2}{\sum_{i=1}^l \alpha_i^2}\leq \max_{1\leq i\leq l} |p(\lambda_i)|\leq \max_{\lambda\in[C_1^{-1},C_2]}|p(\lambda)|.
\end{equation*}
\end{proof}

\end{document}